\documentclass[12pt]{article}
\usepackage{fullpage,indentfirst}
\usepackage{amsfonts,amssymb,stmaryrd}
\usepackage{diagrams}

\def\rk{{\rm rk}}

\def\GB{{\cal B}}
\def\GP{{\cal P}}
\def\GH{{\cal H}}
\def\H{{\bf H}}
\def\B{{\bf B}}
\def\TH{{T^H}}
\def\Chi{{\cal X}}

\def\a{\alpha}
\def\Tau{\mathcal{T}}

\def\QQ{\mathbb{Q}}
\def\PP{\mathbb{P}}
\def\CC{\mathbb{C}}

\def\longto{\longrightarrow}

\newtheorem{theo}{Theorem}
\newtheorem{prop}{Proposition}[section]
\newtheorem{lemma}{Lemma}[section]
\newtheorem{coro}{Corollary}[section]

\def\proof{\noindent{\bf Proof:~}}
\def\carre{\hfill $\square$}
\def\vsp{\vspace{3.5mm}}
\def\definition{\noindent{\bf Definition. }}

\newcounter{paragraf}[section]
\renewcommand{\theparagraph}{\thesection.\arabic{paragraf}}

\newcommand{\parag}
{\refstepcounter{paragraf}
\vsp

{\bf \theparagraph}\hspace{0.2em}--- 
}

\begin{document}


\vsp\vspace*{1cm}

\begin{center}
  {\large\bf
About Knop's action of the Weyl group\\ 
on the set of  orbits of a spherical subgroup\\ 
in the flag manifold}\\

\vspace{3mm}
N. Ressayre
\end{center}

\section{Introduction}

Let $G$ be a complex connected reductive algebraic group.
Let $\GB$ denote the flag variety of $G$.  
Let $H$ be an algebraic subgroup of $G$ which has a finite number of orbits
in $\GB$ ; $H$ is said to be {\it spherical}.
We denote by $\H(\GB)$ the set of the $H$-orbits in $\GB$. 
The closures of these orbits are of importance in representation
theory (see \cite{Wolf}). 
Moreover, the elements of $\H(\GB)$, viewed as orbits of a Borel
subgroup of $G$ in $G/H$ play an important role in the geometry and
topology of the $G$-equivariant embeddings $X$ of $G/H$.

In \cite{Kn:WGH}, F. Knop introduced an action of a monoid 
(constructed from the Weyl group of $G$) on $\H(\GB)$. 
This action is called ``weak order'' and studied by
M. Brion in \cite{Br:GammaGH}.
But, the most spectacular combinatoric structure of the set $\H(\GB)$
was discovered by F. Knop in \cite{Kn:WGH}:
he defined  an action of the Weyl group $W$ of $G$
on $\H(\GB)$.
Actually, the results of F. Knop are stated in a more general context.
The proof of the existence of this action is very indirect and sophisticated.
The aim of this note is to construct natural invariants separating 
the $W$-orbits.
Note that our methods are elementary.

Let us fix a maximal torus $T^H$ of $H$. Denote by $W_H$ the Weyl
group of $T^H$. Let $T$ be a maximal torus of $G$ containing $T^H$ and
let $W$ denote the Weyl group of $T$. 

Let $V\in \H(\GB)$. 
Let $x$ be a point of $V$ whose the orbit by $T^H$ is of minimal dimension.
Denote by $S$ the identity component of the stabilizer of $x$ in
$T^H$. 
The group $W_H$ acts naturally on the set of subtori of $T^H$. 
 The $W_H$-orbit of $S$ is called {\it the type of $V$}. 
It is shown in Section~\ref{sec:deftype} that 
the type of $V$ only depends on $V$ and not on $x$. 

The main result of this note is the following\\

{\bf Theorem}
{\it Two elements of $\H(\GB)$ are in the same $W$-orbit for Knop's
action if and only if they have the same type.}\\ 

In Section~\ref{sec:def}, we recall some useful definitions about a graph 
with vertices the elements of $\H(\GB)$, Knop's action of $W$ on $\H(\GB)$ 
and some classical invariants associated to the elements of $\H(\GB)$.
In Section~\ref{sec:deftype}, we show that the definition of the type 
of an orbit of $H$ is consistent. After, we study the fixed points
of subtori of $H$ in the elements of $\H(\GB)$.
In Section~\ref{sec:actiontype}, we state and prove our main results.
In the following one, we give some consequences of our results and our proofs.

\section{Definitions and notation}
\label{sec:def}

\parag
Let us fix some general notation.
If $\Gamma$ denotes a linear algebraic group, we denote by
$\Gamma^\circ$ its identity component.
If $\Gamma$ acts on an algebraic variety $X$ and $x$ belongs to $X$, 
we denote by $\Gamma_x$ the stabilizer of $x$ and by $\Gamma.x$ 
the orbit of $x$. 
The set of  points of $X$ fixed by $\Gamma$ is denoted by $X^\Gamma$.  
If $S$ is a subgroup of $\Gamma$, we denote by $N_\Gamma(S)$ the
normalizer of $S$ in $\Gamma$ and by $\Gamma^S$ the centralizer of $S$
in $\Gamma$. 

\parag
Let us recall that $G$ is a connected complex reductive group,
$\GB$ its flag variety 
and $H$ a closed subgroup of $G$. 
We assume that $H$ is {\it spherical}; that is, $H$ has a dense 
orbit in $\GB$. In this article,
we are interested in the set $\H(\GB)$ of the orbits of $H$ in
$\GB$. It is shown in  \cite{Br:BGHfini}, \cite{Vi:BGHfini} or
\cite{Kn:WGH} that $\H(\GB)$ is finite. 

We recall the definition of \cite{GammaGH} of a graph $\Gamma(G/H)$
whose vertices are the elements of $\H(\GB)$. The original
construction of $\Gamma(G/H)$ due to M. Brion is very slightly different 
(see \cite{Br:GammaGH}).

Consider the set $\Delta$ of conjugacy classes of minimal non solvable
parabolic subgroups of $G$. 
If $\a$ belongs to $\Delta$, we denote by $\GP_\a$ the $G$-homogeneous
space with isotropy $\a$. 
Then, there exists a unique $G$-equivariant map
$\phi_\a\,:\,\GB\longto\GP_\a$ which is a $\PP^1$-bundle. 

Let $V\in\H(\GB)$ and $\a\in\Delta$. 
We assume that the restriction of $\phi_\a$ to $V$ is finite and 
we denote  its degree by $d(V,\alpha)$. 
Then, there exists a unique open $H$-orbit $V'$ in
$\phi_\a^{-1}(\phi_\a(V))$; in this case, we say that {\it $\a$ raises
$V$ to $V'$}.  
One of the following three cases occurs.
\begin{itemize}
\item Type $U$: $H$ has two orbits in $\phi_\a^{-1}(\phi_\a(V))$
($V$ and $V'$) and $d(V,\a)=1$.
\item Type $T$: $H$ has three orbits in $\phi_\a^{-1}(\phi_\a(V))$
 and $d(V,\a)=1$. 
\item Type $N$: $H$ has two orbits in $\phi_\a^{-1}(\phi_\a(V))$
($V$ and $V'$) and $d(V,\a)=2$.
\end{itemize}

\vsp
\definition
Let $\Gamma(G/H)$ be the oriented graph with vertices the elements of
$\H(\GB)$ and edges labeled by $\Delta$, where $V$ is joined to $V'$
by an edge labeled by $\a$ if $\a$ raises $V$ to $V'$. This edge is
simple (resp. double) if $d(V,\a)=1$ (resp. 2). 
Following the above cases, we say that an edge has {\it type $U$,
  $T$ or $N$}.\\  

One can find examples of graphs $\Gamma(G/H)$ in 
\cite{Br:GammaGH,Spin:these,GammaGH}.

\parag 
Let us fix a Borel subgroup $B$ of $G$, and a maximal torus $T$ of $B$. 
Let $W$ denote the Weyl group of $T$.
We now describe Knop's action of $W$ on the set $\H(\GB)$
(see also \cite{Kn:WGH}). 
Indeed, the action of simple
reflexions easily reads off the graph $\Gamma(G/H)$.  

Every $\a$ in $\Delta$ has a unique representative $P_\a$ which
contains $B$. 
Moreover, there exists a unique $s_\a$ in $W$ such that $Bs_\a B$ is
dense in $P_\a$; and this $s_\a$ is a simple reflexion of $W$. 
The map, $\Delta\longto W,\,\a\longmapsto s_\a$ is a bijection from
$\Delta$ onto the set of simple reflexions of $W$.

Consider the group $\widetilde{W}$ generated by $\{s_\a\ :\ \a\in\Delta\}$
with the relations $s_\a^2=1$.
There is a surjective homomorphism $\widetilde{W}\longto W$.
Let $\Tau$ denote its kernel.

One defines an action of $\widetilde{W}$ on the set $\H(\GB)$ by describing 
the action of the $s_\a$, for any $\a\in\Delta$:

\begin{itemize}
\item Type $U$: $s_\a$ exchanges the two vertices of an edge of type
  $U$ labeled by $\a$.
\item Type $T$: If $\a$ raises $V_1$ and $V_2$ on $V$, then
  $s_\a V_1=V_2$ and $s_\a V=V$.
\item Type $N$: $s_\a$ fixes the two vertices of a double edge
  labeled by $\a$.
\item $s_\a$ fixes all others vertices of $\Gamma(G/H)$.
\end{itemize}

In \cite{Kn:WGH}, F. Knop showed that this action of $\widetilde{W}$
factors through $W$; that is , that $\Tau$ acts trivially on $\H(\GB)$.
The aim of this paper is to describe the orbits of this action
by a natural invariant and to give some consequences.
 
\parag
Denote by $\GH$ the $G$-homogeneous space $G/H$. If $V$ belongs 
to $\H(\GB)$, we set:
$$
V_\GH=\{gH/H\,:\,g^{-1}B/B\in V\}.
$$ 
Then, $V_\GH$ is a $B$-orbit in $\GH$. Moreover, the map 
$V\longmapsto V_\GH$ is a bijection from $\H(\GB)$ onto the set
$\B(\GH)$ of $B$-orbits in $\GH$. 

The {\it character group $\Chi(V_\GH)$} of $V$ (or $V_\GH$) is the set of all
characters of $B$ that arise as weights of eigenvectors of $B$ in the
function field $\CC(V_\GH)$.  Then $\Chi(V_\GH)$ is a free abelian group of
finite rank $\rk(V_\GH)$ (or $\rk(V)$), {\it the rank of $V$}.

\section{The type of an orbit of $H$}
\label{sec:deftype}

\parag
In this section, we  define the type of a $H$-orbit in general 
(not only in $\GB$).  
We start with two technical lemmas.

Let us fix a maximal torus $\TH$ of $H$.
If $V$ is a $H$-homogeneous space, we set: 
$$
\rho_V=\min_{x\in V} \dim (\TH.x).
$$
 
\begin{lemma}
\label{lem:StormaxHx}
Let $V\in\H(\GB)$.
Then, for all $x\in V$, the following are equivalent:
\begin{enumerate}
\item $\dim (\TH.x)=\rho_V$,
\item $(\TH_x)^\circ$ is a maximal torus of $H_x$.
\end{enumerate}
\end{lemma}

\proof
Assume that $\dim (\TH.x)=\rho_V$. Let $S'\supseteq (\TH_x)^\circ$ be a maximal
torus of $H_x$. Then, there exists $h$ in  $H$ such that $h^{-1}S'h$ is
contained is $\TH$. But, $h^{-1}S'h$ fixes $h^{-1}x$. 
Therefore, 
$\dim T^H-\dim T^H_x=\rho_V\leq\dim(T^H.h^{-1}x)\leq\dim T^H-\dim S'$; hence
$\dim S'\leq\dim \TH_x$. 
It follows that $S'=(\TH_x)^\circ$. 

The converse is obvious since $(\TH_x)^\circ$ is always a torus of
$H_x$. 
\carre

\vsp
\begin{lemma}
\label{lem:deftype}
Let $x$ and $y$ belong to $V$ such that $\dim (\TH.x)=\dim
(\TH.y)=\rho_V$.  
Set $S_x=(\TH_x)^\circ$ and $S_y=(\TH_y)^\circ$.   

Then, we have: 
\begin{enumerate} 
\item 
\label{ass:deftype1}
There exists $h$ in $H$ such that $y=h.x$ and 
  $S_y =h S_x h^{-1}$.  
\item \label{ass:deftype2}
There exist $\hat{w}\in  N_H(\TH)$ such that 
$\hat{w}^{-1}S_y\hat{w}=S_x$ and  
$\hat{w}^{-1}.y\in H^{S_x}.x$.
\end{enumerate}   
\end{lemma}  

\proof  
Let $h_1\in H$ such that $y=h_1.x$.  
By Lemma~\ref{lem:StormaxHx}, $h_1^{-1}S_yh_1$ and $S_x$ are
maximal tori of  $H_x=h_1^{-1}H_yh_1$.  
Therefore, (see \cite[21.3]{Hum}) there exists $h_2$ in $H_x$ such that 
 $h_2^{-1}h_1^{-1}S_yh_1h_2=S_x$. Then, $h=h_1h_2$ satisfies 
Assertion~\ref{ass:deftype1}.

Notice that $H^{S_x}=h^{-1}H^{S_y}h$. Then, $\TH$ and $h^{-1}\TH h$ are
maximal tori of $H^{S_x}$; so there exists $g_1$ in $H^{S_x}$ such
that $g_1^{-1}h^{-1}\TH hg_1=\TH$. But, we have:
$g_1^{-1}h^{-1}S_yhg_1=S_x$. 
Then, $\hat{w}=hg_1$ satisfies Assertion~\ref{ass:deftype2}. 
\carre

\vsp
Let $W_H=N_H(\TH)/\TH$ denote the Weyl group of $H$.  
The group $W_H$ acts by conjugacy on the set of subtori of $\TH$.   
Let  $V$ be a $H$-homogeneous space.  
Let us fix $x$ in $V$ such that $\rho_V=\dim(T^H.x)$. 
Then, by Lemma~\ref{lem:deftype}, the orbit $W_H.(T^H_x)^\circ$ does  
not depend on $x$ but only on $V$; we call it {\it the type of $V$}.

\parag
We have:

\begin{prop}
\label{prop:VS} 
Let $S$ belong to the type of $V$. 
Then, we have:
\begin{enumerate}
\item \label{ass:VS1}
$V^S$ is a unique orbit of $N_H(S)$.
\item \label{ass:VS2}
The irreducible components of $V^S$ are orbits of
  $(H^S)^\circ$. 
\end{enumerate}
\end{prop}

\proof
Since $V$ is stable by $H$, $V^S$ is stable by $N_H(S)$. 
Let $x$ and $y$ belong to $V^S$. 
Let $h\in H$  such that $y=h.x$. 
Then, $h^{-1}Sh$ is contained in $H_x$. So by Lemma \ref{lem:StormaxHx}, 
$S$ and $h^{-1}Sh$ are maximal tori of  $H_x$ and hence there exists $h_1$ in
$H_x$ such that $h_1^{-1}h^{-1}Shh_1=S$.
Then,  $y=hh_1.x$ belongs to $N_H(S).x$. 
Assertion~\ref{ass:VS1} is  proved.

By \cite[Corollary 16.3]{Hum}, the identity component of  $N_H(S)$ is
$(H^S)^\circ$. Then, Assertion~\ref{ass:VS2} follows 
from Assertion~\ref{ass:VS1}. 
\carre

\section{The type of an orbit of $H$ in $\GB$}

\parag
\label{par:choixB}
In the previous section, we associated to each $H$-homogeneous space 
$V$ a type and an integer $\rho_V$.
Now, we apply these constructions to the orbits $V$ of $H$ in $\GB$.
First, Proposition~\ref{prop:chiV} below shows that the type of $V$
corresponds to the character group of $V$.
We will deduce that $\rho_V-{\rm rk}(V)$ is independent of $V$.

Let us fix a maximal torus $T$ of $G$ containing $\TH$. 
Let $B$ be a Borel subgroup of $G$ containing $T$.

\begin{prop}
\label{prop:chiV}
Let $V$ be in $\H(\GB)$ and  $S$ be a subtorus of $T$ which belongs 
to the type of $V$. 
Let $w\in W$ such that $V$ intersects the irreducible component
$G^S.wB/B$ of $\GB^S$. 

Then,  $\Chi(V)\otimes\QQ$ is equal to $\Chi(T)^{w^{-1}Sw}\otimes \QQ$.
\end{prop}

\proof 
Let $g\in G$ such that $gB/B$ belongs to $V\cap G^S.wB/B$. 
Consider $y=g^{-1}H/H$. 
By replacing $g$ by an element of $gB$, we may assume that
$\dim(T.y)=\min_{y'\in B.y}\dim(T.y')$. 
But, by Lemma~\ref{lem:StormaxHx} $T_y^\circ$ is a maximal torus of $B_y$. 
Since the unipotent radical of $B_y^\circ$ is contained in $U$, it is
equal to $U_y$. Then, we have: $G_y^\circ=T_y^\circ U_y$.

We have: $\Chi(V)\otimes\QQ=\Chi(B)^{B_y^\circ}\otimes \QQ$. 
Moreover, the restriction map from $\Chi(B_y^\circ)$ to
$\Chi(T_y^\circ)$ is injective. 
Therefore, $\Chi(V)\otimes\QQ=\Chi(T)^{T_y^\circ}\otimes \QQ$. 

Since $B_y=g^{-1}H_xg$, $g^{-1}Sg$ is a maximal torus of $B_y$. 
Therefore, there exists $b\in B_y^\circ$ such that $S=gbT_y^\circ
b^{-1}g^{-1}$.
By replacing $g$ by $gb$ (and keeping $x$ and $y$ unchanged), we may
assume that $b$ is trivial; that is, that 
$S=gT_y^\circ g^{-1}$.

It follows that $T$ and $gTg^{-1}$ are maximal tori of $G^S$. Then,
there exists $s\in G^S$ such that 
$sg$ normalizes $T$. 
Let $w_1$ be the class of $sg$ in the Weyl group of $T$.
Then, $T_y^\circ =w_1^{-1}Sw_1$. 

On the other hand, since $sg\in G^S wB$, there exists $w'$ in the 
Weyl group of $G^S$ such that $w_1=w'.w$. 
Then, $T_y^\circ =w^{-1}Sw$ and the proposition follows.
\carre

\begin{coro}
\label{cor:rkrho}
Let $V$ be an orbit of $H$ in $\GB$. We have:
\begin{enumerate}
\item \label{ass:rkrho1}
$\rho_V-\rk(V)=\rk(G)-\rk (H)$.
\item The rank of $V$ is minimal in $\H(\GB)$ if and only if $V$
  contains points fixed by $\TH$. 
\end{enumerate}
\end{coro}

\proof
The proposition shows that the rank of $V$ is the dimension of $T$
minus the dimension of $S$. On the other hand, $\rho_V$ is the
difference between the rank of $H$ and the dimension of $S$. 
Assertion~\ref{ass:rkrho1} follows. 

Since $\TH$ has fixed points in $\GB$, the rank of $V$ is minimal if
and only if $\rho_V=0$; that is, if and only if $V$ contains points
fixed by $T^H$.
\carre

\parag
Let $V$ be in $\H(\GB)$ and $S$ belong to the type of $V$.
We are now interested in the set $V^S$.
We can make Proposition~\ref{prop:VS} more precise:

\begin{prop}
\label{prop:VSGB}
  \begin{enumerate}
  \item The intersection of $V^S$ and an irreducible component of
  $\GB^S$ is a unique  orbit of $H^S$.
\item If $H$ is connected, the intersection of $V$ and one
  irreducible component of $\GB^S$ is irreducible. 
  \end{enumerate}
\end{prop}

\proof
Let $x$ and $y$ be two points of $V^S$ in the same irreducible component of
  $\GB^S$.
Since the irreducible components of $\GB^S$ are orbits of $G^S$, 
there exists  $g\in G^S$  such that $y=g.x$.
By Assertion $(i)$ there exists $h\in N_H(S)$ such that $y=h.x$.
Then, $g^{-1}h$ belongs to $G_x$ which is a Borel subgroup of $G$ which
contains $S$. Moreover, $g^{-1}h$ normalizes $S$. 
But, by \cite[Proposition 19.4]{Hum}, we have: $N_{G_x}(S)=G_x^S$. 
So, $g^{-1}h$ and $h$ centralize $S$. Assertion $(iii)$ follows. 

If $H$ is connected, Theorem 22.3 of \cite{Hum} shows that $H^S$ is
connected. Now, Assertion $(iv)$ follows from Assertion $(iii)$.
\carre

\parag
We are now interested in the set of irreducible components of $\GB^S$
which intersect $V$. 
By Proposition~\ref{prop:VSGB}, if $H$ is connected, this set is in 
bijection with the set of the irreducible components of $V^S$.

Since the irreducible components of $\GB^S$ are the $G^SwB/B$ 
for $w$ in $W$, we set:
$$
\mathcal{C}(V,S)=\{w\in W\,:\,V\cap G^S wB/B\neq\emptyset\}.
$$
To describe $\mathcal{C}(V,S)$, we need two technical lemmas.

\begin{lemma}
\label{lem:WHS}
Set $N_H(S)G^S=\{hg\,:\,h\in N_H(S){\rm~and~}g\in G^S\}$.

Then, $N_H(S)G^S$ is a closed subgroup of $N_G(S)$ whose 
identity component is $G^S$. 
Moreover, the group $(N_H(S)G^S)/G^S$ is isomorphic to 
$N_H(S)/H^S$ (the Weyl group of $S$ in $H$, denoted by
$W(H,S)$).
\end{lemma}

\proof
Notice that, $N_H(S)$ normalizes $G^S$. Now, one easily checks that
$N_H(S)G^S$ is a subgroup of $G$. 
Moreover,  $N_H(S)G^S$ is clearly contained in $N_G(S)$ and
contains $G^S$. But by \cite[Corollary 16.3]{Hum}, $G^S$ is the
identity component of $N_G(S)$. 
It follows that the index of $G^S$ in $N_H(S)G^S$ is finite. 
Then, $N_H(S)G^S$ is closed in $N_G(S)$ and  its
identity component is $G^S$.  
The last assertion is obvious.
\carre

\vsp

Notice that $T$ is contained in $N_H(S)G^S$. 
Set $W_{N_H(S)G^S}=N_{N_H(S)G^S}(T)/T$. 
Then,  the inclusion of $N_{N_H(S)G^S}(T)$ in $N_G(T)$ induces an  
embedding of $W_{N_H(S)G^S}$ in $W$. 
Let $W_{G^S}$ denote the Weyl group of $(G^S,T)$.

\begin{lemma}
\label{lem:exactseq}
We have an exact sequence:
$$
1\longto W_{G^S}\longto
W_{N_H(S)G^S}\longto
W(H,S)\longto 1.
$$  
\end{lemma}

\proof
Let us start with the exact sequence given by Lemma~\ref{lem:WHS}:
$$
1\longto G^S\longto
N_H(S)G^S\longto
W(H,S)\longto 1.
$$  
By intersecting with $N_{N_H(S)G^S}(T)$, we obtain an exact sequence:
$$
1\longto N_{G^S}(T)\longto
N_{N_H(S)G^S}(T)\longto
W(H,S),
$$  
and it is sufficient to prove that the last map is surjective.
Let $h$ in $N_H(S)$ and $g$ in $G^S$.
Since, $ghT(gh)^{-1}$ and $T$ are maximal tori of $G^S$, there exists
$g'\in G^S$ such that $g'ghT(gh)^{-1}g'^{-1}=T$.
The lemma follows.
\carre\\

If $E$ is a finite set, let $|E|$ denote its cardinality.
Now, we can describe $\mathcal{C}(V,S)$:

\begin{prop}
\label{prop:WVSWHtilde}
\begin{enumerate}
\item The set $\mathcal{C}(V,S)$  is an orbit of $W_{N_H(S)G^S}$ for
its action on $W$ by left multiplication. 
\item If $H$ is connected, $V^S$ has $|W_{N_H(S)G^S}|$ 
irreducible components.
\end{enumerate}
\end{prop}

\proof
Let $\sigma$ be an element of $\mathcal{C}(V,S)$ and let $x$ belong to 
$V\cap G^S\sigma B/B$. 
By Proposition \ref{prop:VS}, $V^S=N_H(S).x$. 
Therefore $G^S.V^S=G^SN_H(S).x=(N_H(S)G^S) \sigma B/B$. 
But $G^SV^S$ is the union of the $G^S.wB/B$ for $w\in \mathcal{C}(V,S)$. 
The first assertion follows. 

By Proposition~\ref{prop:VSGB}, each irreducible component of $V^S$
is the intersection of $V$ and one irreducible component of 
$\GB^S=\coprod_{w\in {}_{W_{G^S}}\backslash W}G^SwB/B$.
Therefore, by the first assertion $V^S$ has 
$\frac{|W_{N_H(S)G^S}|}{|W_{G^S}|}$ irreducible components.
Now, the second assertion follows from Lemma~\ref{lem:exactseq}.
\carre

\parag
Each irreducible component of $\GB^S$ is isomorphic to the flag variety
$\GB_{G^S}$ of $G^S$. 
Moreover, by Proposition~\ref{prop:VSGB}, $V$ intersects any such 
irreducible component in one orbit of $H^S$.
We will now describe the orbits of $H^S$ in $\GB_{G}$ which appear 
in that way.

Let $\tau$ be a $W_H$-orbit of subtori of $\TH$.
Let $\H(\GB)_\tau$ denote the set of $H$-orbits in $\GB$ of type
$\tau$. 

\begin{prop}
\label{prop:HStoH}
Assume that $\H(\GB)_\tau$ is not empty.
Let us fix an element $S$ in $\tau$. Then,
\begin{enumerate}
\item The subgroup $H^S$ of $G^S$ is spherical. 
\item 
The rank of $G^S/H^S$ is equal to the rank of the free abelian 
group $\Chi(T)^S$. 
\item\label{ass:HStoH1}
Let $V\in \H(\GB)_\tau$ and $x\in V^S$.
Then, $\rho_{H^S.x}=\rk(H)-\rk(S)$.
In particular, $\rk(H^S.x)=\rk(G^S/H^S)$.
\item 
\label{ass:HStoH2}
Conversely, let $y$ in $\GB^S$ such that $\rho_{H^S.y}=\rk(H)-\rk(S)$. 
Then, the type of $H.y$ is  $\tau$. 
\end{enumerate} 
\end{prop}

\proof
We first prove Assertions~\ref{ass:HStoH1} and \ref{ass:HStoH2}.
Let $V\in \H(\GB)_\tau$ and $x\in V^S$.\\

Let $y\in H^S.x$. Since $y$ belongs to $V$ and the type of $V$ is $\tau$, we
have $\dim(\TH_y)\leq\dim S$. 
Then, $\rho_{H^S.y}\leq {\rm rk}(H)-{\rm rk}(S)$.

But $\rho_{H^S.x}\geq \rho_{H.x}={\rm rk}(H)-{\rm rk}(S)$.
So $\rho_{H^S.x}={\rm rk}(H)-{\rm rk}(S)$.
This proves Assertion~\ref{ass:HStoH1}.

Set $\Omega=\{y\in G^S.x\ :\ \rho_{H^S.y}\leq {\rm rk}(H)-{\rm rk}(S)\}$.
The set $\Omega$ is open in $G^S.x$ and contains $x$.

Let $y\in\Omega$.
Then, $S$ is a maximal torus of $H^S_y$. 
Let $S_y$ be a maximal torus of $H_y$ containing $S$. 
Then, $S_y$ is contained in $H^S$. Therefore $S=S_y$. 
Then, Lemma \ref{lem:StormaxHx} shows 
that $\rho_{H.y}={\rm rk}(H)-{\rm rk}(S)$.
Therefore, since $(H.y)^S$ is not empty, the type of $H.y$ is
$\tau$.
By Corollary~\ref{cor:rkrho}, this proves Assertion~\ref{ass:HStoH2}.\\

By Proposition~\ref{prop:VSGB}, each orbit of type $\tau$ intersects
$G^S.x$ in a unique orbit of $H^S$. 
Hence, Assertion~\ref{ass:HStoH1} shows that the set of $H^S$-orbit in 
$\Omega$ is finite. So, $H^S$ has a dense orbit in $\Omega$ and in $G^S.x$.
The first assertion follows.
The second one is now a consequence of Assertion~\ref{ass:HStoH1}.
\carre

\section{Knop's action of $W$ on $\H(\GB)$ and orbit  type} 
\label{sec:actiontype}

\parag
Keep the notation as above.
In particular, $\tau$ is a $W_H$-conjugacy class of subtori of $\TH$ such that
$\H(\GB)_\tau$ is not empty and $S$ belongs to $\tau$. 
Set 
${}_{W_{N_H(S)G^S}}\backslash W=\{W_{N_H(S)G^S}w\;:\;w\in W\}$.
By Proposition~\ref{prop:WVSWHtilde}, we can define a map   
$$
\begin{array}{r@{}ccl}
\Theta\,:\,&\H(\GB)_\tau&\longto&{}_{W_{N_H(S)G^S}}\backslash W\\
&V&\longmapsto &\mathcal{C}(V,S).
\end{array}
$$
We consider on ${}_{W_{N_H(S)G^S}}\backslash W$ the action of
the Weyl group $W$ by right multiplication. 

In this section we  show the following 

\begin{theo}
\label{th:ThetaWeq}
The subset $\H(\GB)_\tau$ of $\H(\GB)$ is stable by  Knop's action
of $W$. Moreover, the map $\Theta$ is $W$-equivariant.
\end{theo}

\parag 
Start with    

\begin{lemma}
\label{lem:2alt}
Let $V\in\H(\GB)_\tau$, $x\in V^S$ and $\a\in \Delta$. 
Consider $\phi_\a\,:\,\GB\longto\GP_\a$.
Let $w\in W$ be such that $G^S.x=G^S.wB/B$. 
Then one of the two following cases occurs: 

\begin{tabular}{ll}
\underline{Case 1}:~& $\phi_\a^{-1}(\phi_\a(x))$ is pointwise fixed by $S$.\\
                   & Then, we have $G^Sws_\a B/B=G^SwB/B$.\\
\underline{Case 2}:~&
There exists $y\neq x$ such that  $\phi_\a^{-1}(\phi_\a(x))^S=\{x,y\}$.\\
& Then, $G^S.x\neq G^S.y$ and $G^S.y=G^Sws_\a B/B$.
\end{tabular}
\end{lemma}

\proof
Set $F=\phi_\a^{-1}(\phi_\a(x))$. The variety $F$ is isomorphic to the
projective line $\PP^1$. 
Moreover, $F$ is stable by the action of the torus $S$. 
Then, the  image of $S$ in Aut$(F)\simeq PSL(2)$ is
either trivial or a maximal torus of Aut$(F)$. 
In particular, one of the following cases occurs.

\begin{tabular}{ll}
Case 1:~&$F^S=F$.\\
Case 2:~&There exists $y\neq x$ such that $F^S=\{x,y\}$.
\end{tabular}

In either case, consider the $G^S$-orbit $G^S.\phi_\a(x)$ and the flag variety
$\GB_{G^S}$ of the group $G^S$. 
Since $G^S.\phi_\a(x)$ is the image by
$\phi_\a$ of $G^S.x\simeq \GB_{G^S}$, it is a complete
$G^S$-homogeneous space. Moreover, since $\phi_\a$ is a
$\PP^1$-fibration, we have:  
$\dim(\GB_{G^S})\geq\dim(G^S.\phi_\a(x))\geq\dim(\GB_{G^S})-1$.
Then, two cases occur.  

\begin{tabular}{ll}
Case a:~ &$G^S_{\phi_\a(x)}$ is a non solvable minimal parabolic  
subgroup of  $G^S$\\
& and $G^S.x$ contains $F$.\\
Case b:~ &$G^S_{\phi_\a(x)}$ is a Borel subgroup of $G^S$ 
and $F\cap G^S.x=\{x\}$. 
\end{tabular}

In Case 1, $F$ is contained in the irreducible component of $\GB^S$
which contains $x$; that is in  $G^S.x$. So, Case 1 implies Case a.
In Case 2, we cannot have that $F$ contains  $G^S.x$. 
So, Case 2 implies Case b. In particular, $G^S.x\neq G^S.y$. 

It remains to determine $G^S.ws_\a B/B$ in each case.  
The fiber $\phi_\a^{-1}(\phi_\a(B/B))$ of $\phi_\a$ is the closure
$\overline{Bs_\a B}/B$ of  $Bs_\a B/B$ in $\GB$. 
Let $g\in G^S$ be such that $x=gwB/B$. 
Then, $F=gw\overline{Bs_\a B}/B$. 

In Case 1, $F$ is contained in $G^SwB/B$. In particular, $gws_\a$
belongs to $G^SwB/B$. Therefore, $G^Sws_\a B/B=G^SwB/B$.

In Case 2, we can notice that $gws_\a B/B$ is fixed by $S$ and belongs
to $F$; Therefore, $y=gws_\a B/B$. 
Then, $G^S.y=G^Sws_\a B/B$.
\carre

\parag
{\bf Proof of Theorem \ref{th:ThetaWeq}.}
Let $V\in\H(\GB)_\tau$ and $\a\in\Delta$. 
We will prove that $\mathcal{C}(V,S)s_\a =\mathcal{C}(s_\a V,S)$. 
Let $w\in \mathcal{C}(V,S)$. 
By Proposition~\ref{prop:WVSWHtilde}, it is sufficient to
show that $ws_\a$ belongs to $\mathcal{C}(s_\a V,S)$. 

We fix $x$ in $V^S\cap G^SwB/B$ and we set $F=\phi_\a^{-1}(\phi_\a(x))$. 
Then, one of the following 4 cases occurs. 

\vsp
\noindent\underline{Case 1:} 
$\a$ raises $V$ on $s_\a V$ (type $U$).

Since $V\cap F=\{x\}$, $(s_\a V)^S$ is not empty.
Since $V$ and $s_\a V$ have the same rank, 
Corollary~\ref{cor:rkrho} implies that $S$ belongs to the type of $s_\a V$.

Let us assume that there exists $y\neq x$ such that $F^S=\{x,y\}$. 
Necessarily, $y$ belongs to $s_\a V$. 
But Lemma~\ref{lem:2alt} shows that  $G^S.y=G^Sws_\a B/B$. 
So, $ws_\a$ belongs to $\mathcal{C}(s_\a V,S)$.

If $F^S=F$ then Lemma~\ref{lem:2alt} shows that $F$ is contained in 
$G^SwB/B=G^Sws_\a B/B$. 
Then, since $F$ intersects $s_\a V$, $ws_\a$ belongs 
to $\mathcal{C}(s_\a V,S)$.

\vsp
\noindent\underline{Case 2:}
$\a$ raises $V$ and $s_\a V$ on a third $H$-orbit $V_1$ (type $T$). 

Since $\rk(V_1)=\rk(V)+1$, Corollary~\ref{cor:rkrho} shows that
$V_1^S$ is empty. Then, by Lemma \ref{lem:2alt} there exists $y\neq x$
such that $F^S=\{x,y\}$. 
On the other hand, $F\cap V_1$ is equal to $F$ with two points
removed (type $T$). 
Since, $V_1^S$ is empty it follows that $F\cap V_1=F-\{x,y\}$, 
$F\cap V=\{x\}$ and $F\cap s_\a V=\{y\}$. 
But, Lemma \ref{lem:2alt} shows that $G^S.y=G^Sws_\a B/B$. 
Therefore, $ws_\a$ belongs to $\mathcal{C}(s_\a V,S)$.

\vsp
\noindent\underline{Case 3:}
$\a$ raises $V$ and $s_\a V=V$ (type $N$).

The same proof as in Case 2 shows that $F^S=\{x,y\}=F\cap V$ and
$G^S.y=G^Sws_\a B/B$. It follows that 
$ws_\a$ belongs to $\mathcal{C}(V,S)=\mathcal{C}(s_\a V,S)$.

\vsp
\noindent\underline{Case 4:} $F\cap V$ is open in $F$.

If $F^S=F$ then $V$ is the only $H$-orbit in $\phi_\a^{-1}(V)$ of
maximal rank (type $T$ or $N$). 
Therefore, $s_\a V=V$. Moreover, by Lemma \ref{lem:2alt}, we have  
$G^S.wB/B=G^Sws_\a B/B$; 
therefore, $ws_\a\in \mathcal{C}(V,S)=\mathcal{C}(s_\a V,S)$.

We may assume that there exists $y\neq x$ such that $F^S=\{x,y\}$. 
Then, since $F\cap V$ is open in $F$, stable by $S$ and contains $x$, 
$F\cap V$ is either $F$ or $F-\{y\}$. 
If $F\cap V=F-\{y\}$ then $\a$ raises $s_\a V$ to $V$ by an edge of
type $U$. 
By exchanging $V$ and $s_\a V$ we come back to Case 1. 
Assume that $V$ contains $F$. Since $G^Sws_\a B/B$ intersects $F$, it
intersects $V$. 
Then, $V=s_\a V$ and $ws_\a\in \mathcal{C}(V,S)=\mathcal{C}(s_\a V,S)$.

This completes the proof of Theorem \ref{th:ThetaWeq}.
\carre

\parag
Let $\sigma$ be in $W$ and $\overline{\sigma}$ be its class in 
${}_{W_{N_H(S)G^S}}\backslash W$.
We are now interested in the fiber $\Theta^{-1}(\overline{\sigma})$ 
of $\Theta$.
By definition of $\mathcal{C}(V,S)$, 
$\Theta^{-1}(\overline{\sigma})$ is the set of the orbits $V$ in 
$\H(\GB)_\tau$ which intersects $G^S\sigma B/B$.
Let $\H^S(\GB_{G^S})$ denote the set of the $H^S$-orbits in the flag 
manifold $\GB_{G^S}$ of $G^S$,
and let $\H^S(\GB_{G^S})_{\rm max}$ denote the set of the $H^S$-orbits 
of maximal rank.
By Proposition~\ref{prop:HStoH}, the map
$$
\begin{array}{l@{}ccl}
\eta_\sigma\ :\ &\Theta^{-1}(\overline{\sigma})&\longto&
\H^S(\GB_{G^S})_{\rm max}\\
&V&\longmapsto&V\cap G^S\sigma B/B
\end{array}
$$
is a bijection.

The subgroup $\sigma^{-1}W_{N_H(S)G^S}\sigma$ stabilizes 
$\Theta^{-1}(\overline{\sigma})$.
Moreover, $W_{N_H(S)G^S}$ contains $W_{G^S}$.
Therefore, the group $W_{G^S}$ acts on 
$\Theta^{-1}(\overline{\sigma})$ through the morphism
$W_{G^S}\longto W,\ w\longmapsto\sigma^{-1}w\sigma$.
On the other hand, $W_{G^S}$ acts on $\H^S(\GB_{G^S})_{\rm max}$
by Knop's action.
Is the bijection $\eta_\sigma$ $W_{G^S}$-equivariant ?
The answer is NO in general, but YES for at least one $\sigma$.

\begin{prop}
\label{prop:etaeq}
There exists $\sigma$ such that $\eta_\sigma$ is $W_{G^S}$-equivariant.
\end{prop}

\proof
Actually, the map $\Theta$ depends on the choice of the Borel subgroup
$B$ made in Paragraph~\ref{par:choixB}.
To prove the proposition, it is sufficient to prove that for a good choice 
of $B$, $\eta_1$ is $W_{G^S}$-equivariant.
Let us make such a choice.

Let $P$ be a parabolic subgroup of $G$ with  Levi subgroup  $G^S$.
Let $B$ be a Borel subgroup of $G$ such that $T\subset B\subset P$. 

Notice that $B^S=B\cap G^S$ is a Borel subgroup of $G^S$.
Denote by $\Delta^S$ the set of conjugacy classes of minimal non
solvable parabolic subgroups of $G^S$.  
Let $\a\in\Delta^S$ and $\GP_\a^S$ denote the $G^S$-homogeneous space
with isotropy $\a$.
If $P_\a^S$ is a minimal parabolic subgroup of $G^S$ containing $B^S$ 
corresponding to $\a$, then  $P_\a^S.B$ is a minimal parabolic subgroup 
of $G$.
Moreover, $P_\a^S=(P_\a^S.B)\cap G^S$.
Therefore, we obtain an immersion (from now on implicit) of $\Delta^S$
in $\Delta$. 
In particular $P_\a=P_\a^SB$.
Consider the following commutative diagram $\mathcal{D}$:

\begin{diagram}
\GB_{G^S}\simeq &G^SB/B&\rInto^{\rm inclusion}&\GB\\
&\dTo&&\dTo_{\phi_{\a_i}}\\
\GP_\a^S\simeq&\GP_{\a_i}^S&\rInto&\GP_{\a_i}
\end{diagram}

The restriction of $\phi_\a$ to $G^SB/B$ is obviously the unique 
$G^S$-equivariant map $\phi_{\a,S}$ from $\GB_{G^S}$ onto $\GP_{\a,S}$.

Let $x\in G^S B/B$ such that $H^S.x$ belongs to $\H^S(\GB_{G^S})_{\rm max}$. 
It remains to prove the following

\begin{center}
  \underline{Claim}:~$G^SB/B\cap (s_\a.Hx)=s_\a(H^Sx)$.
\end{center}

Since Diagram $\mathcal{D}$ is commutative, we have 
\begin{eqnarray}
  \label{eq:phi}
  \phi_\a^{-1}(\phi_\a(x))=
\phi_{\a,S}^{-1}(\phi_{\a,S}(x));
\end{eqnarray}

we denote by $F$ this subvarity of $\GB$.
Moreover, since the rank of $H^Sx$ is maximal in $\H^S(\GB_{G^S})$,
Proposition~\ref{prop:HStoH} shows

\begin{eqnarray}
  \label{eq:GSB}
G^SB/B\cap Hx=H^Sx.  
\end{eqnarray}

Four cases can occur:

\begin{tabular}{ll}
Case 1:~&$\a$ raises $H^Sx$ in $\Gamma(G^S/H^S)$.\\
Case 2:& $\a$ raises an orbit $H^Sy$ of $\H^S(\GB_{G^S})_{\rm max}$
on $H^Sx$.\\
Case 3:&$\a$ raises an orbit of $\H^S(\GB_{G^S})$ on $H^Sx$ 
by an edge of type $T$ or $N$.\\ 
Case 4:& $H^S x=\phi_{\a,S}^{-1}(\phi_{\a,S}(H^Sx))$.
\end{tabular}

In Case 1, $F\cap H^Sx=\{x\}$ and $H^S_x$, and hence $H_x$, 
acts transitively on $F-\{x\}$.
Moreover, by Equality~\ref{eq:GSB}, $F\cap Hx=\{x\}$. 
Therefore, $\a$ raises $Hx$ by an edge of type $U$ in $\Gamma(G/H)$.
The claim follows.

In Case 2, $H^Sy$ is in Case 1. 
The claim follows.

In Case 3, $F\cap H^Sx=F\cap Hx$ is equal to $F$ with two points removed.
Therefore, $\a$ raises an orbit of $\H(\GB)$ on $Hx$ by an edge of type $T$ 
or $N$, and $s_\a.Hx=Hx$.

In Case 4, $F$ is contained in $H^Sx$ and hence in $Hx$.
As a consequence, $Hx=\phi_\a^{-1}(\phi_\a(Hx))$ and
$s_\a.Hx=Hx$.
This completes the proof of the proposition.
\carre\\

Here, comes our main result.

\begin{theo}
\label{th}
Two elements of $\H(\GB)$ are in the same $W$-orbit for  Knop's
action if and only if they have the same type.  
\end{theo}

\proof
By Theorem~\ref{th:ThetaWeq}, it is sufficient to prove that one (or any)  
fiber of $\Theta$ is an orbit of $W_{N_H(S)G^S}$.
Then, by Proposition~\ref{prop:etaeq}, it is sufficient to prove the 
theorem for the orbits of maximal rank.
Let $V_0$ be such an orbit.
There exist a sequence $\a_1,\cdots,\a_k$ in $\Delta$ and a sequence
$V_0,\,V_1,\cdots,V_k$ of $H$-orbits such that $\a_i$ raises $V_{i-1}$ 
on $V_i$ for all $i=1,\cdots,k$, and $V_k$ is the open $H$-orbit in $\GB$.
Since the rank of $V_0$ is maximal, all the orbits $V_i$ have the same 
rank and the edges joining these orbits are of type $U$.
Therefore, we have $(s_{\a_k}\cdots s_{\a_1}).V_0=V_k$.
The theorem is proved.
\carre\\

\section{Some consequences}

\parag
Theorem \ref{th} has a nice corollary about the character groups of
the elements of $\B(\GH)$:

\begin{coro}
Let $V$ and $V'$ in $\H(\GB)$. 
Then,
$\Chi(V)=\Chi(V')$ if and only if $\Chi(V)\otimes\QQ=\Chi(V')\otimes\QQ$.
\end{coro}

\proof
Let us fix $T\subset B$.
By identifying $\Chi(B)$ with $\Chi(T)$, we obtain an action of $W$ on
$\Chi(B)$. 
Assume that $\Chi(V)\otimes\QQ=\Chi(V')\otimes\QQ$. 
By Proposition \ref{prop:chiV}, the orbits $V$ and $V'$ have the same
type. 
Then, by Theorem \ref{th}, there exists $w$ in $W$ such that $V=wV'$. 
Then, by \cite[Theorem 4.3]{Kn:WGH}, $\Chi(V)=w.\Chi(V')$.  
Now, $\Chi(V)\otimes\QQ=\Chi(V')\otimes\QQ$ implies
$\Chi(V)=\Chi(V')$. 
\carre 

\parag
We can also apply Theorem 2 to the description of the isotropy subgroups
of the action of $H$ in $\GB$.

\begin{coro}
\label{cor:stab}
  Let $x$ and $y$ be in $\GB$ such that $Hx$ and $Hy$ have the same type.
Then, $(H_x/H_x^\circ)$ and $(H_y/H_y^\circ)$ are isomorphic.
\end{coro}

\proof
Set $V=Hx$ and $V'=Hy$.
Let $\a\in Delta$. 
Since $W$ is generated by the simple reflections, by Theorem~\ref{th}
it is sufficient to prove the corollary for $V'=s_\a.V\neq V$.
Two cases occur:

\begin{itemize}
\item Type $T$: $V$ and $V'$ are raised on a third orbit $V''$.
\item Type $U$: $\a$ raises $V$ on $V'$ (up to re-indexing).
\end{itemize}

In the first case, the restrictions of $\phi_\a$ to $V$ and $V'$ are 
isomorphisms onto $\phi_\a(V'')$.
The corollary follows.

Assume that $\a$ raises $V$ on $V'=s_\a.V$.
By replacing $y$ by another point of $Hy$, we may assume that 
$\phi_\a(x)=\phi_\a(y)$.
Since the restriction of $\phi_\a$ to $V$ is an isomorphism onto $\phi_\a(V)$
and $\phi_\a(V')=\phi_\a(V)$, $H_y$ is contained in $H_x$.
This inclusion induces a morphism 
$\psi\ :\ H_y/H_y^\circ\longto H_x/H_x^\circ$.
But, $H_x/H_y$ is isomorphic to $\mathbb{A}^1$ and hence irreducible. 
We deduce that $\psi$ is surjective.

It remains to show that $H_y\cap H_x^\circ=H_y^\circ$ to prove that $\psi$ 
is injective.
Obviously, $H_y^\circ\subset (H_y\cap H_x^\circ)$; and we can define a
morphism $H_x^\circ/H_y^\circ\longto H_x^\circ/(H_y\cap H_x^\circ)$.
Since $H_x^\circ/(H_y\cap H_x^\circ)$ is isomorphic to $\mathbb{A}^1$, it
is simply connected and $H_y\cap H_x^\circ=H_y^\circ$.
\carre

\parag
We are going to apply Theorem~\ref{th} to the $H$-orbits in $\GB$
of minimal rank.
We keep notation as above.
In particular, $\H(\GB)_{\{\TH\}}$ is the set of the orbits of $H$
in $\GB$ of minimal rank.

\begin{prop}
\label{prop:min}
  We assume that $H$ is connected.
Then, we have:
\begin{enumerate}
\item \label{ass:min1}
The group $H^\TH/\TH$ is a maximal unipotent subgroup of $G^\TH/\TH$.
\item \label{ass:min2}
The stabilizers in $W$ (for Knop's action) of the elements of 
$\H(\GB)_{\{\TH\}}$ are isomorphic to the Weyl group $W_H$ of $H$.
\item Let $V$ be in $\H(\GB)_{\{\TH\}}$.
The stabilizers in $H$ of the points of $V$ are connected.
\end{enumerate}
\end{prop}

\proof
Since $\TH$ is maximal in $H$, $H^\TH/\TH$ is unipotent.
But it is  a spherical subgroup of $G^\TH/\TH$.
Assertion~\ref{ass:min1} follows.

We claim that the cardinality of the set $\H(\GB)_{\{\TH\}}$ is 
$\frac{|W|}{|W_H|}$.
By Proposition~\ref{prop:WVSWHtilde}, we have to prove that the set of 
irreducible components of the $V^\TH$ for $V\in\H(\GB)_{\{\TH\}}$ 
has the same cardinality as $W$.
But, by Proposition~\ref{prop:HStoH}, this set is in natural bijection with
the set of orbits of $H^\TH$ in $\GB^\TH$.
Moreover, by Assertion~\ref{ass:min1}, $H^\TH$ has $|W_{G^\TH}|$ orbits
in each one of the  $\frac{|W|}{|W_{G^\TH}|}$ irreducible components 
of $\GB^\TH$. The claim follows.

By Proposition~\ref{prop:etaeq}, we may assume that $\eta_1$ is 
$W_{G^\TH}$-equivariant to prove Assertion~\ref{ass:min2}.
Let $V$ be in $\H(\GB)_{\{\TH\}}$ such that $\Theta(V)=\overline{1}$.
We have to prove that the stabilizer $W_V$ of $V$ in $W$ is 
isomorphic to $W_H$.
Since $\Theta$ is $W$-equivariant, $W_V$ is contained in 
$W_{N_H(\TH)G^\TH}$ and by Lemma~\ref{lem:exactseq} maps on $W_H$.
Moreover, the claim shows that $|W_V|=|W_H|$.
So, by  Lemma~\ref{lem:exactseq} it is sufficient to prove that 
$W_V\cap W_{G^S}$ is trivial.
By Proposition~\ref{prop:etaeq}, this is a consequence 
of Assertion~\ref{ass:min1}.

By Corollary~\ref{cor:stab}, it is sufficient to prove the last assertion 
for a closed orbit $V$ of $H$ in $\GB$.
Let $x$ be in $V$.
Since $V$ is closed in $\GB$, it is projective.
So, $H_x$ is a parabolic subgroup of $H$. 
In particular, it is connected. 
\carre

\vspace{4mm}
\noindent
{\bf Acknowledgment:} I am grateful to S. Pin for numerous and very
useful discussions.

\bibliographystyle{smfalpha}
\bibliography{biblio}

\begin{center}
  -\hspace{1em}$\diamondsuit$\hspace{1em}-
\end{center}

\vspace{5mm}
\begin{flushleft}
Nicolas Ressayre\\
Universit{\'e} Montpellier II\\
D{\'e}partement de Math{\'e}matiques\\
Case courrier 051-Place Eug{\`e}ne Bataillon\\
34095 Montpellier Cedex 5\\
France\\
e-mail:~{\tt ressayre@math.univ-montp2.fr}  
\end{flushleft}
\end{document}